\documentclass[12pt]{article}
\usepackage{amssymb}
\usepackage{amsfonts}
\usepackage{amsmath}
\usepackage{graphicx}
\usepackage{geometry}
\usepackage[nodisplayskipstretch,doublespacing]{setspace}
\usepackage{amscd}

\newtheorem{theorem}{\sc{Theorem}}[section]
\newtheorem{corollary}[theorem]{\sc{Corollary}}
\newtheorem{definition}[theorem]{\sc{Definition}}
\newtheorem{lemma}[theorem]{\sc{Lemma}}
\newtheorem{example}[theorem]{\sc{Example}}
\newtheorem{remark}[theorem]{\sc{Remark}}
\newenvironment{proof}[1][Proof]{\noindent\textbf{#1.} }{\ \rule{0.5em}{0.5em}}
\input{tcilatex}
\begin{document}

\title{Good Gradings of Generalized Incidence Rings}
\author{Kenneth L. Price \\
%EndAName
University of Wisconsin Oshkosh}
\date{\today}
\maketitle

\begin{abstract}
This inquiry is based on both the construction of generalized incidence
rings due to Gene Abrams and the construction of good group gradings of
incidence algebras due to Molli Jones. We provide conditions for a
generalized incidence ring to be graded isomorphic to a subring of an
incidence ring over a preorder. We also extend Jones's construction to good
group gradings for incidence algebras over preorders with crosscuts of
length one or two.
\end{abstract}

\section{Overview}

Unless otherwise stated we use multiplicative notation for all semigroup,
monoid, or group operations and the identity is denoted by 1.

Suppose $G$ is a semigroup and $S$ is a ring which does not necessarily
contain a multiplicative identity. We say $S$ is a $G$-graded ring if there
is a direct sum $S=\bigoplus_{a\in G}S_{a}$, as a group under the addition
of $S$, such that $S_{b}S_{c}\subseteq S_{bc}$ for all $b,c\in G$. The
subgroups $S_{a}$, $a\in G$, are called the homogeneous components, the
elements of $\cup _{a\in G}S_{a}$ are called the homogeneous elements, and
every element is a sum of finitely many homogeneous elements. We let $%
\partial s$ be the unique element of $G$ such that $s\in S_{\partial s}$ for
any nonzero homogeneous $s\in S$. The support of $S$ is the set $\limfunc{%
Supp}_{G}S=\left\{ a\in G:S_{a}\neq 0\right\} $. The grading is called
finite if $\limfunc{Supp}_{G}S$ is a finite set.

An important type of matrix algebra grading is a good grading (for example,
see references \cite{BZ}, \cite{DIN}, and \cite{MJ}). This definition
extends easily to incidence algebras (see \cite{MJ}, \cite{MillSpieg}, and 
\cite{ParmSpieg}). In section \ref{gen inc ring section} we state the
definition of balanced relation introduced by Abrams (see \cite{Ab}) and go
over the construction of generalized incidence rings. Good semigroup
gradings of generalized incidence rings are defined in section \ref%
{incidence rings section}. Theorems \ref{Good Grading Theorem} and \ref{All
Good Theorem} are fundamental for our constructions since they categorize
good gradings of generalized incidence rings in terms of homomorphisms from
the relations to the semigroup. Theorem \ref{Hasse Diagram Lemma} shows how
to construct good semigroup gradings of incidence algebras over minimally
connected partial orders.

Suppose $S=\bigoplus_{a\in G}S_{a}$ and $T=\bigoplus_{a\in G}T_{a}$ are $G$%
-graded rings. A homomorphism of $G$-graded rings is a ring homomorphism $%
h:S\rightarrow T$ such that $h\left( S_{a}\right) \subseteq T_{a}$ for all $%
a\in \limfunc{Supp}_{G}S$. An isomorphism which is a homomorphism of $G$%
-graded rings is called an isomorphism of $G$-graded rings. In the case of
matrix algebras there are gradings which are not good gradings but are
isomorphic to good gradings (see \cite[Example 1.3]{DIN}). Isomorphic
gradings for good group gradings of incidence algebras over partial orders
have been studied by Miller and Spiegel (see \cite{MillSpieg}).

In section \ref{Stable Relations Section} we state the definitions of
compression maps and stable relations (see \cite{Price}). If $G$ is a
cancellative monoid then Theorem \ref{compression theorem} shows compression
maps provide a correspondence between good $G$-gradings. Stable relations
are used to describe a class of generalized incidence rings which are
isomorphic to subrings of incidence rings over preorders by Theorem \ref%
{Stable to Transitive Thm}. The isomorphism, which is described in Lemma \ref%
{preserving induced lemma}, is an isomorphism of $G$-graded rings.

Good group gradings are considered in section \ref{Group Gradings Section}.
The main result of this section is Theorem \ref{extend Jones theorem}, which
extends \cite[Theorem 4]{MJ} to good group gradings for incidence algebras
over preorders with crosscuts of length one or two. The conclusion of our
theorem is modified to account for finite gradings. We show our result for
preorders is related to generalized incidence rings in Corollary \ref%
{cross-cut corollary} and Example \ref{stable to partial order example}. We
finish with example \ref{infinite support example}, which describes a
partial order with a minimal element whose incidence algebra does not have
the free-extension property. Unfortunately this is a counterexample to \cite[%
Theorem 4]{MJ}.

\section{Generalized Incidence Rings\label{gen inc ring section}}

\begin{figure}[th]
\begin{center}
\includegraphics[width=3.84in,height=0.88in]{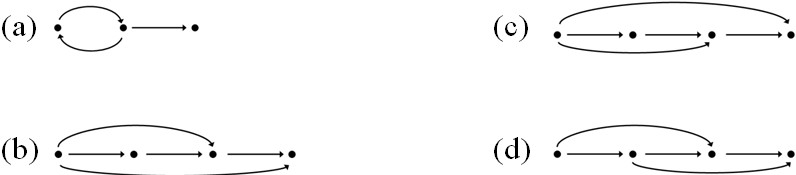}
\end{center}
\caption{Directed graphs determining reflexive relations.}
\label{FourDigraphs}
\end{figure}

By definition a relation $\rho $ on a set $X$ is a subset of $X\times X$. We
adopt the usual convention of writing $x\rho y$ if $x,y\in X$ satisfy $%
\left( x,y\right) \in \rho $. The notation $x\rho y$ is often shorter and
more convenient, but the notation $\left( x,y\right) \in \rho $ will be used
when it is helpful to describe the relation as a set of ordered pairs. The
directed graphs shown in Figure \ref{FourDigraphs} represent reflexive
relations. We omit loops in all diagrams so that the arrows match up with
elements of the off-diagonal subset of relations.

The construction of generalized incidence rings does not require the
relation to be a preorder. To define multiplication we will assume $X$ is
locally finite, that is, every interval of $X$ is a finite set. (Recall an
interval\textit{\ }is a subset of the form $\left[ x,y\right] =\left\{ z\in
X:x\rho z\text{ and }z\rho y\right\} $ for some $x,y\in X$.)

Throughout the rest of this section $\rho $ is a locally finite relation on
a set $X$ and $R$ is an associative ring with unity. Let $I\left( X,\rho
,R\right) $\ denote the set of all functions $f:X\times X\rightarrow R$ such
that $f\left( x,y\right) \neq 0$ implies $x\rho y$. Componentwise operations
determine an $R$-module such that given $r\in R$ and $f,g\in I\left( X,\rho
,R\right) $ the functions $rf,f+g\in I\left( X,\rho ,R\right) $ satisfy $%
\left( rf\right) \left( x,y\right) =r\left( f\left( x,y\right) \right) $ and 
$\left( f+g\right) \left( x,y\right) =f\left( x,y\right) +g\left( x,y\right) 
$ for all $x,y\in X$.

The \emph{balance property} is satisfied by $w,x,y,z\in X$ if any of the
relations $w\rho x$, $x\rho y$, $y\rho z$, and $w\rho z$ do not hold, or all
four relations hold and $w\rho y$ if and only if $x\rho z$. The relation is 
\emph{balanced }if it is reflexive and the balance property is satisfied by
all $w,x,y,z\in X$. The relations determined by directed graphs (a), (b),
and (c) in Figure \ref{FourDigraphs} are not balanced because they are all
missing arrows. Note that a reflexive and transitive relation is balanced
but there are balanced relations, such as the one determined by (d) in
Figure \ref{FourDigraphs}, which are not transitive.

In case $\rho $ is balanced we can combine the proof of \cite[Proposition 1.2%
]{Ab} with the assumption that $X$ is locally finite to construct a ring
multiplication with identity $e:X\times X\rightarrow R$ such that $e\left(
x,x\right) =1$ and $e\left( x,y\right) =0$ for all $x,y\in X$ such that $%
x\neq y$. In this case we call $I\left( X,\rho ,R\right) $ the \emph{%
generalized incidence ring of }$X$\emph{\ with coefficients in }$R$. If $%
\rho $ is a partial order then $I\left( X,\rho ,R\right) $ is the incidence
ring over $R$.

The multiplication is defined so that the product of $f,g\in I\left( X,\rho
,R\right) $ is the function given by equation \ref{function mult} for all $%
x,y\in X$. This is called \emph{convolution}. 
\begin{equation}
\left( fg\right) \left( x,y\right) =\left\{ 
\begin{tabular}{ll}
$\dsum\limits_{z\in \left[ x,y\right] }f\left( x,z\right) g\left( z,y\right) 
$ & if $x\rho y$ in $X$ \\ 
$0$ & otherwise%
\end{tabular}%
\right.  \label{function mult}
\end{equation}%
Obviously if the sum in equation \ref{function mult} is nonzero then $x\rho
y $ in $X$. In this case $\left[ x,y\right] $ is finite so convolution\
uniquely determines an element of $R$. Consider $f,g,h\in I\left( X,\rho
,R\right) $ and $w,z\in X$ with $w\rho z$. A nonzero term of $\left( \left(
fg\right) h\right) \left( w,z\right) $ is determined by $y\in \left[ w,z%
\right] $ and $x\in \left[ w,y\right] $. Thus there are relations $w\rho x$, 
$x\rho y$, $y\rho z$, $w\rho z$, and $w\rho y$. Since the relation is
balanced these are equivalent to $w\rho x$, $x\rho y$, $y\rho z$, $w\rho z$,
and $x\rho z$, or $x\in \left[ w,z\right] $ and $y\in \left[ x,z\right] $.
This gives an identical term of $\left( f\left( gh\right) \right) \left(
w,z\right) $. Since the nonzero terms match up over the relations we have $%
\left( fg\right) h=f\left( gh\right) $.

For each $x,y\in X$ such that $x\rho y$ there exists $e_{xy}\in I\left(
X,\rho ,R\right) $ such that $e_{xy}\left( i,j\right) $ is given by equation %
\ref{standard unit eq} for all $i,j\in X$. 
\begin{equation}
e_{xy}\left( i,j\right) =\left\{ 
\begin{tabular}{ll}
$1$ & if $x=i$ and $y=j$ \\ 
$0$ & otherwise%
\end{tabular}%
\right.  \label{standard unit eq}
\end{equation}%
If $w\in X$ then $w\rho w$ since a balanced relation is reflexive. \ It is
easy to show $\left( e_{ww}\right) ^{2}=e_{ww}$ for all $w\in X$ directly
from the definition in \ref{function mult}. Equations \ref{standard unit
mult} and \ref{balanced mult} also hold for all $f\in I\left( X,\rho
,R\right) $ and $x,y,z,w\in X$ such that $x\rho y$ and $z\rho w$. 
\begin{equation}
f\left( x,y\right) e_{xy}=e_{xx}fe_{yy}  \label{standard unit mult}
\end{equation}%
\begin{equation}
e_{xy}e_{zw}=\left\{ 
\begin{tabular}{ll}
$e_{x,w}$ & if $y=z$ and $x\rho w$ in $X$ \\ 
$0$ & otherwise%
\end{tabular}%
\right.  \label{balanced mult}
\end{equation}

In any situation where we refer to a generalized incidence ring we mean an
associative ring with unity formed on the $R$-module of functions $I\left(
X,\rho ,R\right) $ where $R$ is a ring with unity and $\rho $ is a locally
finite balanced relation on $X$. The operation of $R$ on $I\left( X,\rho
,R\right) $ does not play a significant role in our investigation. We
reserve the term \textit{generalized incidence algebra} for $I\left( X,\rho
,R\right) $ where $R$ is a commutative ring with unity and $\rho $ is a
locally finite balanced relation on $X$. If, additionally, $\rho $ is a
partial order then $I\left( X,\rho ,R\right) $ is the usual incidence
algebra over $R$ (see \cite{SOD}).

\section{Good Gradings\label{incidence rings section}}

\begin{definition}
\label{setup def}Assume $G$ is a semigroup and $\rho $ is a relation on $X$.

\begin{enumerate}
\item Set $\limfunc{Trans}\left( X\right) =\left\{ \left( x,y,z\right)
:x\rho y\text{, }y\rho z\text{, }x\rho z\text{, and }x,y,z\in X\right\} $. A 
\emph{transitive triple in }$X$\emph{\ }is an ordered triple in $\limfunc{%
Trans}\left( X\right) $.

\item We say $\Phi :\rho \rightarrow G$ is a\textit{\ }\emph{homomorphism}
if $\Phi \left( x,y\right) \Phi \left( y,z\right) =\Phi \left( x,z\right) $
holds for any \ $x,y,z\in X$ such that $\left( x,y,z\right) \in \limfunc{%
Trans}\left( X\right) $.

\item If $I\left( X,\rho ,R\right) $ is a $G$-graded generalized incidence
ring then the grading is \emph{good} if $e_{xy}$ is homogeneous for all $%
x,y\in X$ such that $x\rho y$.
\end{enumerate}
\end{definition}

\begin{theorem}
\label{Good Grading Theorem}Assume $G$ is a semigroup, $I\left( X,\rho
,R\right) $ is a generalized incidence ring, and $\Phi :\rho \rightarrow G$
is a homomorphism. Let $S_{a}$ be given by equation \ref{component eq} for
each $a\in G$. Then $I\left( X,\rho ,R\right) =\bigoplus_{a\in G}S_{a}$ is a 
$G$-graded ring if and only if $\func{Im}\Phi $ is finite. 
\begin{equation}
S_{a}=\left\{ f\in I\left( X,\rho ,R\right) :f\left( r\right) \neq 0\text{
implies }\Phi \left( r\right) =a\text{ for all }r\in \rho \right\}
\label{component eq}
\end{equation}
\end{theorem}

%TCIMACRO{\TeXButton{Proof}{\proof}}%
%BeginExpansion
\proof%
%EndExpansion
It is easy to see $S_{a}$ is an $R$-submodule for all $a\in G$ and $%
S_{a}\cap S_{b}=\left\{ 0\right\} $ if $b\in G$ and $b\neq a$. We show $%
S_{a}S_{b}\subseteq S_{ab}$ for all $a,b\in G$. Suppose $f\in S_{a}$, $g\in
S_{b}$ and $\left( fg\right) \left( x,y\right) \neq 0$ for some $x,y\in X$
with $x\rho y$. By equation \ref{function mult} there exists $z\in \left[ x,y%
\right] $ such that $x\rho z$, $z\rho y$, and $f\left( x,z\right) g\left(
z,y\right) \neq 0$. Thus $f\left( x,z\right) \neq 0$ and $g\left( z,y\right)
\neq 0$ which implies $\Phi \left( x,z\right) =a$ and $\Phi \left(
z,y\right) =b$. Moreover, $ab=\Phi \left( x,y\right) $ since $\left(
x,z,y\right) $ is a transitive triple and $\Phi $ is a homomorphism. This
proves $fg\in S_{ab}$ as desired.

To complete the proof we show $\func{Im}\Phi $ is finite if and only if $%
I\left( X,\rho ,R\right) =\bigoplus_{a\in G}S_{a}$. First assume $\func{Im}%
\Phi $ is a finite subset of $G$. Then there is a positive integer $m$ and $%
a_{1},\ldots ,a_{m}\in G$ such that $\func{Im}\Phi =\left\{ a_{1},\ldots
,a_{m}\right\} $. We must prove an arbitrarily chosen $f\in S$ is a sum of
finitely many homogeneous elements.

For each $i=1,\ldots ,m$ let $f_{i}\in I\left( X,\rho ,R\right) $ be the the
function satisfying equation \ref{f comp eq} for all $\left( x,y\right) \in
X $.%
\begin{equation}
f_{i}\left( x,y\right) =\left\{ 
\begin{tabular}{ll}
$f\left( x,y\right) $ & if $x\rho y$ in $X$ and $\Phi \left( x,y\right)
=a_{i}$ \\ 
$0$ & otherwise%
\end{tabular}%
\right.  \label{f comp eq}
\end{equation}%
By construction $f_{i}\in S_{a_{i}}$ for all $i\leq m$. It is easy to prove $%
f$ is the sum of finitely many homogeneous elements, $f_{1},\ldots ,f_{m}$,
as desired.

To prove the other direction assume $S=\bigoplus_{a\in G}S_{a}$. We choose $%
h\in I\left( X,\rho ,R\right) $ so that for all $x,y\in X$ we have $h\left(
x,y\right) =1$ if $x\rho y$ and otherwise $h\left( x,y\right) =0$. If $%
S=\bigoplus_{a\in G}S_{a}$ then we may write $h$ as the sum of homogeneous
elements. Thus there is a positive integer $p$ and distinct $b_{1},\ldots
,b_{p}\in G$ such that $h=h_{1}+\cdots +h_{p}$ where $\partial h_{i}=b_{i}$
for each $i=1,\ldots ,p$. We will show $\func{Im}\Phi \subseteq \left\{
b_{1},\ldots ,b_{p}\right\} $.

For all $x,y\in X$ if $x\rho y$ we have $h\left( x,y\right) =1$ hence $%
h_{1}\left( x,y\right) +\cdots +h_{p}\left( x,y\right) =1$. Then $%
h_{i}\left( x,y\right) \neq 0$ for some $i\leq p$ so, by construction, we
have $\Phi \left( x,y\right) =b_{i}$. Since $x,y\in X$ with $x\rho y$ were
arbitrarily chosen we can conclude $\func{Im}\Phi \subseteq \left\{
b_{1},\ldots ,b_{p}\right\} $. This gives the desired result, $\func{Im}\Phi 
$ is finite.%
%TCIMACRO{\TeXButton{End Proof}{\endproof}}%
%BeginExpansion
\endproof%
%EndExpansion

\begin{definition}
We say \textit{a }\emph{grading is induced by a homomorphism }$\Phi $ if it
can be constructed in the setting of Theorem \ref{Good Grading Theorem} with 
$\func{Im}\Phi $ finite.
\end{definition}

\begin{theorem}
\label{All Good Theorem}Let $S=I\left( X,\rho ,R\right) $ be a generalized
incidence ring and let $G$ be a monoid. Suppose $S=\bigoplus_{a\in G}S_{a}$
is a $G$-graded ring.

\begin{enumerate}
\item If the grading is good then there is a homomorphism $\Phi :\rho
\rightarrow G$ given by $\Phi \left( x,y\right) =\partial e_{x,y}$ for all $%
x,y\in X$ with $x\rho y$.

\item If the grading is good and $e_{xx}\in S_{1}$ for all $x\in X$ then the
grading is induced by $\Phi $. Moreover, the grading is finite.
\end{enumerate}
\end{theorem}

%TCIMACRO{\TeXButton{Proof}{\proof}}%
%BeginExpansion
\proof%
%EndExpansion
Part 1 can be proved directly from equation \ref{balanced mult}. To prove
part 2 we let $a\in \limfunc{Supp}_{G}S$ be arbitrarily chosen. We will show 
$S_{a}$ is given by equation \ref{component eq} so the grading is induced by 
$\Phi $. Then the grading is finite since $\limfunc{Supp}_{G}S=\func{Im}\Phi 
$.

Suppose $f\in I\left( X,\rho ,R\right) $ and for all $x,y\in X$ with $x\rho
y $ if $f\left( x,y\right) \neq 0$ then $\Phi \left( x,y\right) =a$. We must
prove such an $f$ is contained in $S_{a}$. Since $S$ is graded there is a\
positive integer $m$ and nonzero homogeneous $f_{1},\ldots ,f_{m}\in S$ such
that $f=f_{1}+\cdots +f_{m}$ and $\partial f_{i}\neq \partial f_{j}$ if $%
i\neq j$. Fix $i\leq m$. Since $f_{i}$ is nonzero there exist $x,y\in X$
such that $x\rho y$ and $f_{i}\left( x,y\right) \neq 0$. By equation \ref%
{standard unit mult} we have $f_{j}\left( x,y\right)
e_{xy}=e_{xx}f_{j}e_{yy} $ for $j=1,\ldots ,m$. If $f_{j}\left( x,y\right)
\neq 0$ then $\partial e_{xy}=\partial f_{j}$ since $\partial
e_{xx}=\partial e_{yy}=1$. When $j=i$ we conclude $\Phi \left( x,y\right)
=\partial f_{i}$ since $\Phi \left( x,y\right) =\partial e_{xy}$ by
construction and $f_{i}\left( x,y\right) \neq 0$ by our choice of $x,y\in X$%
. But if $j\neq i$ and $f_{j}\left( x,y\right) \neq 0$ then $\partial
e_{xy}=\partial f_{j}$ so $\partial f_{i}=\partial f_{j}$, which is a
contradiction. We are left with $f_{j}\left( x,y\right) =0$ for all $j\neq i$%
. Therefore $f\left( x,y\right) =f_{i}\left( x,y\right) \neq 0$ and $\Phi
\left( x,y\right) =a$ by the assumption on $f$. We already proved $\Phi
\left( x,y\right) =\partial f_{i}$ so $\partial f_{i}=a$. Since $i$ was
arbitrarily chosen we have $\partial f_{i}=a$ for all $i=1,\ldots ,m$.\
Since $\partial f_{1},\ldots ,\partial f_{m}$ are distinct, we must have $%
m=1 $. Thus $f=f_{1}\in S_{a}$, as desired.

Now suppose $g\in S_{a}$ and $g\left( w,z\right) \neq 0$ for some $w,z\in X$
such that $w\rho z$. Equation \ref{standard unit mult} becomes $g\left(
w,z\right) e_{wz}=e_{ww}ge_{zz}$ and we have $\partial e_{wz}=\partial g$
since $\partial e_{ww}=\partial e_{zz}=1$ by assumption. Moreover $%
a=\partial g$ and $\partial e_{wz}=\Phi \left( w,z\right) $ so $\Phi \left(
w,z\right) =a$. 
%TCIMACRO{\TeXButton{End Proof}{\endproof}}%
%BeginExpansion
\endproof%
%EndExpansion

\begin{remark}
Equation \ref{balanced mult} gives $\left( e_{xx}\right) ^{2}=e_{xx}$ and so 
$\Phi \left( x,x\right) ^{2}=\Phi \left( x,x\right) $ for all $x\in X$. The
condition $e_{xx}\in S_{1}$ for all $x\in X$ stated in part 2 of Theorem \ref%
{All Good Theorem} may not hold for all monoids. But if $G$ is a
cancellative monoid then $\left( e_{xx}\right) ^{2}=e_{xx}$ implies $\Phi
\left( x,x\right) =1$ for all $x\in X$ and the condition in part 2 is
satisfied automatically.
\end{remark}

\begin{definition}
Let $G$ be a semigroup and let $\rho $ be a relation on $X$.

\begin{enumerate}
\item A subset $\beta $ of $\rho $ is a $G$\emph{-extendible set}\textit{\
for }$\rho $ if for every function $\phi :\beta \rightarrow G$ there exists
a homomorphism $\Phi :\rho \rightarrow G$ such that $\Phi |_{\beta }=\phi $.
If we may choose $\Phi $ so that $\func{Im}\Phi $ is finite then we say $%
\phi \ $is \emph{grading admissible}.

\item A subset $\gamma $ of $\rho $ is a $G$\emph{-essential set}\textit{\
for }$\rho $ if for all homomorphisms $\Phi _{1},\Phi _{2}:\rho \rightarrow
G $ such that $\Phi _{2}\neq \Phi _{1}$ there exists $c\in \gamma $ such
that $\Phi _{2}\left( c\right) \neq \Phi _{1}\left( c\right) $.

\item A $G$-extendible and $G$-essential subset is called a $G$\emph{%
-grading set}\textit{\ for }$\rho $.
\end{enumerate}
\end{definition}

\begin{remark}
$\sigma $ is a $G$-grading set for $\rho $ if and only if for every function 
$\phi :\sigma \rightarrow G$ there is a unique homomorphism $\Phi :\rho
\rightarrow G$ such that $\Phi |_{\sigma }=\phi $. If $\rho $ is balanced
and $\rho $ contains a $G$-grading set then all good $G$-gradings of a
generalized incidence ring $S=I\left( X,\rho ,R\right) $ such that $%
e_{xx}\in S_{1}$ for all $x\in X$ are uniquely determined by grading
admissible functions from the $G$-grading set to $G$.
\end{remark}

We finish this section with a result on partial orders.\ Recall $\left(
X,\leq \right) $ is minimally connected if $\left( X,\leq \right) $ is a
connected, locally finite partial-order and $\left[ x,y\right] $ is either
empty or a chain for all $x,y\in X$. The Hasse diagram of $X$ is the
directed graph $H$ with vertex set $X$ and arrow set $\left\{ \left(
a,b\right) :a,b\in X\text{ and }b\text{ covers }a\right\} $.

\begin{theorem}
\label{Hasse Diagram Lemma}Assume $G$ is a semigroup, $R$ is an associative
ring with unity, and $\left( X,\leq \right) $ is a minimally-connected
partial order. Then the arrow set of the Hasse diagram of $\left( X,\leq
\right) $ is a $G$-grading set for $\leq $.
\end{theorem}

%TCIMACRO{\TeXButton{Proof}{\proof}}%
%BeginExpansion
\proof%
%EndExpansion
Let $\phi :\sigma \rightarrow G$ be given, where $\sigma $ is the set of
arrows in the Hasse diagram for $\left( X,\leq \right) $. We use $\phi $ to
define a function $\Phi :\leq \rightarrow G$. Let $x,y\in X$ such that $%
x\leq y$ be given. Since $\left( X,\leq \right) $ is a minimally connected
there is a unique chain $\left\{ x_{1},\ldots ,x_{m}\right\} $ in $\sigma $
such that $x=x_{1}\leq x_{2}\leq \ldots \leq x_{m}=y$ and $\left(
x_{i},x_{i+1}\right) \in \sigma $ for each $i<m$. We set $\Phi \left(
x,y\right) =\phi \left( x_{1},x_{2}\right) \cdots \phi \left(
x_{m-1},x_{m}\right) $. A straightforward check proves $\Phi $ is a
homomorphism such that $\Phi |_{\sigma }=\phi $. Thus $\sigma $ is a $G$%
-extendible set for $\leq $.

Suppose $\Phi _{1},\Phi _{2}$ are homomorphism and $\Phi _{1}|_{\sigma
}=\Phi _{2}|_{\sigma }$. Give $x,y\in X$ such that $x\rho y$ there is a
unique chain $\left\{ x_{1},\ldots ,x_{m}\right\} $ in $\sigma $ such that $%
x=x_{1}\leq x_{2}\leq \ldots \leq x_{m}=y$. The homomorphism property gives $%
\Phi _{i}\left( x,y\right) =\phi \left( x_{1},x_{2}\right) \cdots \phi
\left( x_{m-1},x_{m}\right) $ for $i=1$ or $i=2$. Thus $\Phi _{1}=\Phi _{2}$
which proves $\sigma $ is a $G$-essential set for $\rho $. Therefore $\sigma 
$ is a $G$-grading set for $\rho $.%
%TCIMACRO{\TeXButton{End Proof}{\endproof}}%
%BeginExpansion
\endproof%
%EndExpansion

\section{Compression Maps and Stable Relations\label{Stable Relations
Section}}

We fix the notation $\delta \left( \rho \right) =\left\{ \left( x,x\right)
:x\in X\right\} $ for the diagonal subset of a relation $\rho $ on a set $X$%
. The off-diagonal set $\rho ^{\ast }=\rho \backslash \delta \left( \rho
\right) $ is an anti-reflexive relation on $X$.

\begin{definition}
\label{Compression Definition}Suppose $\rho _{1}$ and $\rho _{2}$ are
relations on $X_{1}$ and $X_{2}$, respectively. A function $\theta
:X_{2}\rightarrow X_{1}$ is called a \emph{compression map }if 1, 2, and 3
are satisfied. In this case we say $\rho _{1}$ is a \emph{compression} of $%
\rho _{2}$.

\begin{enumerate}
\item $\theta $ is surjective and order-preserving.

\item For all $a_{1},a_{2},a_{3}\in X_{1}$ if $\left(
a_{1},a_{2},a_{3}\right) \in \limfunc{Trans}\left( X_{1}\right) $ then there
exist $x_{1},x_{2},x_{3}\in X_{2}$ such that $\left(
x_{1},x_{2},x_{3}\right) \in \limfunc{Trans}\left( X_{2}\right) $ and $%
\theta \left( x_{i}\right) =a_{i}$ for $i=1,2,3$.

\item There is a bijection $\theta ^{\ast }:\rho _{2}^{\ast }\rightarrow
\rho _{1}^{\ast }$ given by $\theta ^{\ast }\left( x,y\right) =\left( \theta
\left( x\right) ,\theta \left( y\right) \right) $ for all $x,y\in X_{2}$
with $x\rho _{2}^{\ast }y$.
\end{enumerate}
\end{definition}

Compression maps were introduced in \cite{Price}. Example \ref{stable to
partial order example} uses a compression map to construct a $G$-grading set
for a non-transitive relation.

\begin{lemma}
\label{preserving induced lemma}Suppose $I\left( X_{1},\rho _{1},R\right) $
and $I\left( X_{2},\rho _{2},R\right) $ are generalized incidence rings and $%
\theta :X_{2}\rightarrow X_{1}\ $is a compression map.

\begin{enumerate}
\item If $G$ is a semigroup and $\Phi _{1}:\rho _{1}\rightarrow G$ is a
homomorphism then a homomorphism $\Phi _{2}:\rho _{2}\rightarrow G$ is given
by $\Phi _{2}\left( x,y\right) =\Phi _{1}\left( \theta \left( x\right)
,\theta \left( y\right) \right) $ for all $x,y\in X$ with $x\rho _{2}y$.

\item If $S=I\left( X_{1},\rho _{1},R\right) $ and $T=I\left( X_{2},\rho
_{2},R\right) $ have gradings induced by $\Phi _{1}$ and $\Phi _{2}$,
respectively, then there is an injective homomorphism of $G$-graded rings $%
h:S\rightarrow T$ such that equation \ref{graded homo eq} holds for all $%
f\in S$ and all $x,y\in X_{2}$.%
\begin{equation}
\left( h\left( f\right) \right) \left( x,y\right) =\left\{ 
\begin{tabular}{ll}
$f\left( \theta \left( x\right) ,\theta \left( y\right) \right) $ & if $%
x\rho _{2}y$ in $X_{2}$ \\ 
$0$ & otherwise%
\end{tabular}%
\right.  \label{graded homo eq}
\end{equation}
\end{enumerate}
\end{lemma}

%TCIMACRO{\TeXButton{Proof}{\proof}}%
%BeginExpansion
\proof%
%EndExpansion
(1) Part 1 is trivial. (2) A routine check shows $h$ is an $R$-module
homomorphism. In the definition of multiplication in equation \ref{function
mult}, it is easy to see $h$ is a ring homomorphism if $\theta \left( \left[
x,y\right] \right) =\left[ \theta \left( x\right) ,\theta \left( y\right) %
\right] $ for all $x,y\in X_{2}$ such that $x\rho _{2}y$. For any $c\in %
\left[ \theta \left( x\right) ,\theta \left( y\right) \right] $ there is a
transitive triple $\left( \theta \left( x\right) ,c,\theta \left( y\right)
\right) \in \limfunc{Trans}\left( X_{1}\right) $ so by part 2 of definition %
\ref{Compression Definition} there exists $u,v,w\in X_{2}$ such that $\left(
\upsilon ,v,w\right) \in \limfunc{Trans}\left( X_{2}\right) $, $\theta
\left( u\right) =\theta \left( x\right) $, $\theta \left( v\right) =c$, and $%
\theta \left( w\right) =\theta \left( y\right) $. We find $u=x$ and $w=y$
since $\theta ^{\ast }\left( x,y\right) =\theta ^{\ast }\left( u,v\right) $
and $\theta ^{\ast }$ is bijective. Therefore $v\in \left[ x,y\right] $ and $%
c\in \theta \left( \left[ x,y\right] \right) $. Since $c$ was arbitrarily
chosen we have $\theta \left( \left[ x,y\right] \right) \subseteq \theta
\left( \left[ x,y\right] \right) $. Moreover $\theta \left( \left[ x,y\right]
\right) \subseteq \left[ \theta \left( x\right) ,\theta \left( y\right) %
\right] $ since $\theta $ is order-preserving. Therefore $\theta \left( %
\left[ x,y\right] \right) =\left[ \theta \left( x\right) ,\theta \left(
y\right) \right] $ as desired.

Given $f\in S\backslash \left\{ 0\right\} $ there exist $a_{1},b_{1}\in
X_{1} $ such that $f\left( a_{1},b_{1}\right) \neq 0$. If $a_{1}\neq b_{1}$
then there exist $x_{1},y_{1}\in X_{2}$ such that $x_{1}\rho _{2}^{\ast
}y_{1}$ and $\theta ^{\ast }\left( x_{1},y_{1}\right) =\left(
a_{1},b_{1}\right) \in \rho _{2}^{\ast }$. If $a_{1}=b_{1}$ then there
exists $x_{1}\in X_{2}$ such that $a_{1}=\theta \left( x_{1}\right) $ and we
set $y_{1}=x_{1}$. In either case $x_{1}\rho _{2}y_{1}$ and $h\left(
f\right) \left( x_{1},y_{1}\right) =f\left( a_{1},b_{1}\right) $. Therefore $%
h\left( f\right) \left( x_{1},y_{1}\right) \neq 0$ and $h\left( f\right) $
is nonzero. This proves $h $ is injective.

Let $c\in \limfunc{Supp}_{G}S$ and $g\in S_{c}\backslash \left\{ 0\right\} $
be arbitrarily chosen. Since $h\left( g\right) \in T\backslash \left\{
0\right\} $ there exist $x_{2},y_{2}\in X_{2}$ such that $x_{2}\rho
_{2}y_{2} $ and $\left( h\left( g\right) \right) \left( x_{2},y_{2}\right)
\neq 0$. If we set $a_{2}=\theta \left( x_{2}\right) $ and $b_{2}=\theta
\left( y_{2}\right) $ then $a_{2}\rho _{1}b_{2}$ since $\theta $ preserves
order. Moreover $\left( h\left( g\right) \right) \left( x_{2},y_{2}\right)
=g\left( a_{2},b_{2}\right) $ and $\Phi _{1}\left( a_{2},b_{2}\right) =c$
since the grading on $S$ is induced by $\Phi _{1}$. We have $\Phi _{2}\left(
x_{2},y_{2}\right) =\Phi _{1}\left( \theta \left( x_{2}\right) ,\theta
\left( y_{2}\right) \right) $ so $\Phi _{2}\left( x_{2},y_{2}\right) =c$
hence $h\left( g\right) \in T_{c}$ as desired. This proves $h$ preserves the
grading.%
%TCIMACRO{\TeXButton{End Proof}{\endproof}}%
%BeginExpansion
\endproof%
%EndExpansion

\begin{theorem}
\label{compression theorem}Suppose $\rho _{1}$ and $\rho _{2}$ are reflexive
relations on $X_{1}$ and $X_{2}$, respectively, $\theta \ $is a compression
map of $X_{2}$ onto $X_{1}$, and $G$ is a cancellative monoid. There is a $G$%
-\textit{grading set of }$\rho _{1}$\textit{\ }if and only if there is a $G$%
-grading set of $\rho _{2}$.
\end{theorem}

%TCIMACRO{\TeXButton{Proof}{\proof}}%
%BeginExpansion
\proof%
%EndExpansion
We assume there is a $G$-grading set \textit{of }$\rho _{2}$\textit{\ }and
prove there is a $G$-grading set of $\rho _{1}$. The reverse implication can
be proved using a similar argument. Suppose $\sigma _{2}$ is a $G$-grading
set of $\rho _{2}$ over $G$. Then $\sigma _{2}\subseteq \rho _{2}^{\ast }$
since $G$ is cancellative. We set $\sigma _{1}=\left\{ \theta ^{\ast }\left(
a,b\right) :a,b\in X_{1}\text{ and }\left( a,b\right) \in \sigma
_{2}\right\} $.

First we prove $\sigma _{1}$ is a $G$-essential set for $\rho _{1}$. Suppose 
$\Phi _{1},\Phi _{1}^{\prime }:\rho _{1}\rightarrow G$ are homomorphisms
such that $\Phi _{1}\neq \Phi _{1}^{\prime }$. Then there exist $%
a_{1},a_{2}\in X_{1}$ such that $a_{1}\rho _{1}^{\ast }a_{2}$ and $\Phi
_{1}\left( a_{1},a_{2}\right) \neq \Phi _{1}^{\prime }\left(
a_{1},a_{2}\right) $. We have $\left( a_{1},a_{2}\right) =\theta ^{\ast
}\left( x_{1},x_{2}\right) $ for some $x_{1},x_{2}\in X_{2}$ since $\theta
^{\ast }$ is bijective. Let $\Phi _{2}:\rho _{2}\rightarrow G$ be the
homomorphism defined as in part 1 of Lemma \ref{preserving induced lemma}.
Let $\Phi _{2}^{\prime }:\rho _{2}\rightarrow G$ be defined similarly, using 
$\Phi _{1}^{\prime }$ in place of $\Phi _{1}$. We have $\Phi _{2}\left(
x_{1},x_{2}\right) =\Phi _{1}\left( a_{1},a_{2}\right) $ and $\Phi
_{2}^{\prime }\left( x_{1},x_{2}\right) =\Phi _{1}^{\prime }\left(
a_{1},a_{2}\right) $ so $\Phi _{2}\left( x_{1},x_{2}\right) \neq \Phi
_{2}^{\prime }\left( x_{1},x_{2}\right) $. Thus there exist $y,z\in X_{1}$
such that $\left( y,z\right) \in \sigma _{2}$ and $\Phi _{2}\left(
y,z\right) \neq \Phi _{2}^{\prime }\left( y,z\right) $ since $\sigma _{1}$
is a $G$-grading set for $\varrho _{1}$. We have $\theta ^{\ast }\left(
y,z\right) \in \sigma _{1}$, $\Phi _{1}\left( \theta ^{\ast }\left(
y,z\right) \right) =\Phi _{2}\left( y,z\right) $, and $\Phi _{1}^{\prime
}\left( \theta ^{\ast }\left( y,z\right) \right) =\Phi _{2}^{\prime }\left(
y,z\right) $ by construction. Therefore $\Phi _{1}\left( \theta ^{\ast
}\left( y,z\right) \right) \neq \Phi _{1}^{\prime }\left( \theta ^{\ast
}\left( y,z\right) \right) $ and this shows $\sigma _{1}$ is a $G$-essential
set for $\rho _{1}$.

Next we prove $\sigma _{1}$ is a $G$-extendible set for $\rho _{1}$. Given $%
\psi _{1}:\sigma _{1}\rightarrow G$ we let $\psi _{2}:\sigma _{2}\rightarrow
G$ be the function given by $\psi _{2}=\psi _{1}\circ \theta ^{\ast }$.
There is also a homomorphism $\Psi _{2}:\rho _{2}\rightarrow G$ such that $%
\Psi _{2}|_{\sigma _{2}}=\psi _{2}$ since $\sigma _{2}$ is a $G$-extendible
set for $\rho _{2}$. There is a homomorphism $\Psi _{1}:\rho _{1}\rightarrow
G$ such that for all $a,b\in X_{1}$ such that $a\rho _{1}b$ we have $\Psi
_{1}\left( a,b\right) =\Psi _{2}\left( \left( \theta ^{\ast }\right)
^{-1}\left( a,b\right) \right) $ if $a\neq b$ and$\Psi _{1}\left( a,b\right)
=1$ if $a=b$. Given $a,b\in X_{1}$ such that $\left( a,b\right) \in \sigma
_{1}$ there exist $x,y\in X_{1}$ such that $a=\theta \left( x\right) $, $%
b=\theta \left( y\right) $, and $\left( x,y\right) \in \sigma _{2}$. We have 
$\theta ^{\ast }\left( x,y\right) =\left( a,b\right) $ and $\Psi _{2}\left(
x,y\right) =\Psi _{1}\left( a,b\right) $ by construction, $\Psi _{2}\left(
x,y\right) =\psi _{2}\left( x,y\right) $ since $\left( x,y\right) \in \sigma
_{2}$, and $\psi _{2}\left( x,y\right) =\psi _{1}\left( a,b\right) $ since $%
\psi _{2}=\psi _{1}\circ \theta ^{\ast }$. This shows $\Psi _{1}|_{\sigma
_{1}}=\psi _{1}$ and $\sigma _{1}$ is a $G$-extendible set for $\rho _{1}$.%
%TCIMACRO{\TeXButton{End Proof}{\endproof}}%
%BeginExpansion
\endproof%
%EndExpansion

\begin{definition}
\label{Stable Definition}Let $\rho $ be a reflexive relation on a set $X$.

\begin{enumerate}
\item $\left( X,\rho \right) $ is \emph{stable }if\ $\rho $ is balanced and
if the relations $a\rho b$, $a\rho c$, $b\rho c$, $b\rho d$, and $c\rho d$
imply $a\rho d$ for all distinct $a,b,c,d\in X$.

\item An element $x\in X$ is a \emph{clasp} if there exist $w,y\in
X\backslash \left\{ x\right\} $ such that $w\rho x$, $x\rho y$, and $\left(
w,y\right) \notin \rho $.

\item $x\in X$ is a \emph{locked clasp} if there exist $u,v,w,y\in
X\backslash \left\{ x\right\} $ such that $\left( w,y\right) \notin \rho $
and $\left( u,x,y\right) ,\left( u,x,v\right) ,\left( w,x,v\right) \in 
\limfunc{Trans}\left( X\right) $ .

\item An \emph{unlocked clasp} is a clasp which is not locked.
\end{enumerate}
\end{definition}

It is easy to see a preorder is stable. The balanced relation determined by
(d) in Figure \ref{FourDigraphs} is not stable. Neither a balanced relation
which is not stable nor a stable relation which contains a locked clasp can
be the compression of a preorder by \cite[Theorem 2.4 and Lemma 3.4]{Price}.

\begin{theorem}
\label{Stable to Transitive Thm}Assume $\rho $ is a stable relation on a
finite set $X$ and every clasp in $X$ is unlocked.

\begin{enumerate}
\item There is a preorder $\leq $ on a finite set $Y$ and a compression map $%
\theta :Y\rightarrow X$.

\item There is an injective ring homomorphism $h:I\left( X,\rho ,R\right)
\rightarrow I\left( Y,\leq ,R\right) $.

\item If $G$ is a semigroup and $I\left( X,\rho ,R\right) $ has a grading
induced by $\Phi :\rho \rightarrow G$ then there is a good $G$-grading of $%
I\left( Y,\leq ,R\right) $ such that $h:I\left( X,\rho ,R\right) \rightarrow
I\left( Y,\leq ,R\right) $ is an injective homomorphism of $G$-graded rings.
\end{enumerate}
\end{theorem}

%TCIMACRO{\TeXButton{Proof}{\proof}}%
%BeginExpansion
\proof%
%EndExpansion
Part 1 is \cite[Theorem 3.5]{Price}. Part 2 is a special case of part 3 when 
$G$=$\left\{ 1\right\} $. Part 3 follows from part 1 and Lemma \ref%
{preserving induced lemma}.%
%TCIMACRO{\TeXButton{End Proof}{\endproof}}%
%BeginExpansion
\endproof%
%EndExpansion

\section{Group Gradings\label{Group Gradings Section}}

If $\left( X,\rho \right) $ is a preorder then a subset $A$ of $X$ is a
cross-cut if $A$ is an antichain in $X$, for all $x\in X$ there exists $a\in
A$ such that $x\rho a$ or $a\rho x$, and if $C$ is a chain in $X$ then $C$
can be extended to a chain $C^{^{\prime }}$ such that $A\cap C^{\prime }$ is
nonempty. The length of a cross-cut $A$ of $X$ is $\left\vert A\right\vert $%
. For example, the minimal elements of $X$ form a cross-cut.

For a partially ordered with a cross-cut of length one or two there is a
subset $\sigma $ of $\rho $ such that $\sigma $ is a $G$-grading set for any
group $G$ (see \cite[Theorem 4]{MJ}). However if the shortest cross-cut of a
partially ordered set has length three or more then there may not be a $G$%
-grading set for a group $G$ (see \cite[Example 6]{MJ}).

\begin{theorem}
\label{extend Jones theorem}Assume $\left( X,\rho \right) $ is a preorder
and $X$ has a cross-cut of length one or two. Then there is a subset $\sigma 
$ of $\rho $ such that $\sigma $ is a $G$-grading set for any group $G$.
\end{theorem}

%TCIMACRO{\TeXButton{Proof}{\proof}}%
%BeginExpansion
\proof%
%EndExpansion
If $r,s\in X$ satisfy $r\rho s$ and $s\rho r$ then $r$ and $s$ are said to
be paired in $X$. This defines an equivalence relation on $X$. We set $%
\tilde{X}=\left\{ \left[ x\right] :x\in X\right\} $ where $\left[ x\right]
=\left\{ y\in X:x\text{ and }y\text{ are paired}\right\} $ for all $x\in X$.
There is a partial order $\tilde{\rho}$ on $\tilde{X}$ such that $\left[ x%
\right] \tilde{\rho}\left[ y\right] $ if and only if $x\rho y$ for all $%
x,y\in X$. It is easy to see $\tilde{X}$ also has a cross-cut of length one
or two. A $G$-grading set for $\left( \tilde{X},\tilde{\rho}\right) $ is
constructed in the proof of \cite[Theorem 4]{MJ}. We use it to construct a $%
G $-grading set for $\left( X,\rho \right) $. Fix $P\subseteq X$ such that $%
\tilde{X}=\left\{ \left[ y\right] :y\in P\right\} $ and $\left[ y_{1}\right]
\neq \left[ y_{2}\right] $ for all $y_{1},y_{2}\in P$ such that $y_{1}\neq
y_{2}$. Let $\beta $ be a subset of $\rho $ such that $\beta \subseteq
P\times P$ and $\tilde{\beta}=\left\{ \left( \left[ a\right] ,\left[ b\right]
\right) :\left( a,b\right) \in \beta \right\} $ is a $G$-grading set for $%
\left( \tilde{X},\tilde{\rho}\right) $. Set $\sigma =\beta \cup \gamma $
where $\gamma =\bigcup\limits_{p\in P}\left\{ \left( p,x\right) :x\in \left[
p\right] \backslash \left\{ p\right\} \right\} $.

Let $\phi :\sigma \rightarrow G$ be given. Define $\psi :\tilde{\beta}%
\rightarrow G$ by $\psi \left( \left[ a\right] ,\left[ b\right] \right)
=\phi \left( a,b\right) $ for all $a,b\in P$ such that $\left( a,b\right)
\in \beta $. Then there exists a homomorphism $\Psi :\tilde{\rho}\rightarrow
G$ such that $\Psi |_{\tilde{\beta}}=\psi $. We extend $\phi $ to $\bar{\phi}%
:\sigma \cup \left\{ \left( p,p\right) :p\in P\right\} \rightarrow G$ so
that $\bar{\phi}|_{\sigma }=\phi $ and $\bar{\phi}\left( p,p\right) =1$ for
all $p\in P$. Suppose $x_{1},x_{2}\in X$ such that $x_{1}\rho x_{2}$ are
given. There are uniquely determined $p_{1},p_{2}\in P$ such that $\left[
x_{1}\right] =\left[ p_{1}\right] $ and $\left[ x_{2}\right] =\left[ p_{2}%
\right] $. We set $\Phi \left( x_{1},x_{2}\right) =\bar{\phi}\left(
p_{1},x_{1}\right) ^{-1}\Psi \left( \left[ p_{1}\right] ,\left[ p_{2}\right]
\right) \bar{\phi}\left( p_{2},x_{2}\right) $. A routine check proves the
function $\Phi :\rho \rightarrow G$ is a homomorphism.

Suppose $x_{1},x_{2}\in X$ such that $\left( x_{1},x_{2}\right) \in \sigma $
are given. If $\left( x_{1},x_{2}\right) \in \beta $ then $x_{1},x_{2}\in P$%
, $\Phi \left( x_{1},x_{2}\right) =\bar{\phi}\left( x_{1},x_{1}\right)
^{-1}\Psi \left( \left[ x_{1}\right] ,\left[ x_{2}\right] \right) \bar{\phi}%
\left( x_{2},x_{2}\right) $, $\Psi \left( \left[ x_{1}\right] ,\left[ x_{2}%
\right] \right) =\psi \left( \left[ x_{1}\right] ,\left[ x_{2}\right]
\right) $, and $\psi \left( \left[ x_{1}\right] ,\left[ x_{2}\right] \right)
=\phi \left( x_{1},x_{2}\right) $ by construction. This gives $\Phi \left(
x_{1},x_{2}\right) =\phi \left( x_{1},x_{2}\right) $. If $\left(
x_{1},x_{2}\right) \in \gamma $ then $\Phi \left( x_{1},x_{2}\right) =\bar{%
\phi}\left( x_{1},x_{1}\right) ^{-1}\Psi \left( \left[ x_{1}\right] ,\left[
x_{1}\right] \right) \bar{\phi}\left( x_{1},x_{2}\right) $ since $x_{1}\in P$
and $\left[ x_{1}\right] =\left[ x_{2}\right] $. We have $\bar{\phi}\left(
x_{1},x_{2}\right) =\phi \left( x_{1},x_{2}\right) $ since $\left(
x_{1},x_{2}\right) \in \gamma $ and $\Phi \left( x_{1},x_{2}\right) =\phi
\left( x_{1},x_{2}\right) $ follows easily. We have shown $\Phi |_{\sigma
}=\phi $ so $\sigma $ is a $G$-extendible set for $\rho $.

Suppose $\Upsilon _{1},\Upsilon _{2}:\rho \rightarrow G$ are homomorphism
such that $\Upsilon _{1}|_{\sigma }=\Upsilon _{2}|_{\sigma }$. For $i=1,2$
define $\tilde{\Upsilon}_{i}:\tilde{\rho}\rightarrow G$ by $\tilde{\Upsilon}%
_{i}\left( \left[ a\right] ,\left[ b\right] \right) =\Upsilon _{i}\left(
a,b\right) $ for all $a,b\in P$ such that $\left( a,b\right) \in \rho $. A
routine check proves $\tilde{\Upsilon}_{1},\tilde{\Upsilon}_{2}$ are
homomorphisms such that $\tilde{\Upsilon}_{1}|_{\tilde{\beta}}=\tilde{%
\Upsilon}_{2}|_{\tilde{\beta}}$. Thus $\tilde{\Upsilon}_{1}=\tilde{\Upsilon}%
_{2}$ since $\tilde{\beta}$ is a $G$-grading set for $\tilde{X}$. This
proves $\Upsilon _{1}\left( a,b\right) =\Upsilon _{2}\left( a,b\right) $ for
all $a,b\in P$ such that $\left( a,b\right) \in \rho $.

Suppose $x_{1},x_{2}\in X$ such that $x_{1}\rho x_{2}$ are given. There are
uniquely determined $p_{1},p_{2}\in P$ such that $\left[ x_{1}\right] =\left[
p_{1}\right] $ and $\left[ x_{2}\right] =\left[ p_{2}\right] $. Transitivity
and the homomorphism property gives $\Upsilon _{1}\left( x_{1},x_{2}\right)
=\Upsilon _{1}\left( p_{1},x_{1}\right) ^{-1}\Upsilon _{1}\left(
p_{1},p_{2}\right) \Upsilon _{1}\left( p_{2},x_{2}\right) $. By the result
of the previous paragraph $\Upsilon _{1}\left( p_{1},p_{2}\right) =\Upsilon
_{2}\left( p_{1},p_{2}\right) $. Moreover $\Upsilon _{1}\left(
p_{2},x_{2}\right) =\Upsilon _{2}\left( p_{2},x_{2}\right) $ and $\Upsilon
_{1}\left( p_{1},x_{1}\right) =\Upsilon _{2}\left( p_{1},x_{1}\right) $
since $\Upsilon _{1}|_{\sigma }=\Upsilon _{2}|_{\sigma }$ and $\left(
p_{1},x_{1}\right) ,\left( p_{2},x_{2}\right) \in \gamma $. Substitution
gives $\Upsilon _{1}\left( x_{1},x_{2}\right) =\Upsilon _{2}\left(
p_{1},x_{1}\right) ^{-1}\Upsilon _{2}\left( p_{1},p_{2}\right) \Upsilon
_{2}\left( p_{2},x_{2}\right) $ and this reduces to $\Upsilon _{1}\left(
x_{1},x_{2}\right) =\Upsilon _{2}\left( x_{1},x_{2}\right) $. Since $%
x_{1},x_{2}\in X$ such that $\left( x_{1},x_{2}\right) \in \rho $ were
arbitrarily chosen we can conclude $\Upsilon _{1}=\Upsilon _{2}$ and $\sigma 
$ is a $G$-essential set for $\rho $. Therefore $\sigma $ is a $G$-grading
set for $\rho $. 
%TCIMACRO{\TeXButton{End Proof}{\endproof}}%
%BeginExpansion
\endproof%
%EndExpansion

\begin{corollary}
\label{cross-cut corollary}Assume $\left( X,\rho \right) $ is the
compression of a preorder $\left( Y,\leq \right) $. If $Y$ has a cross-cut
of length one or two then there is a subset $\sigma $ of $\rho $ such that $%
\sigma $ is a $G$-grading set\ for $\rho $ for any group $G$.
\end{corollary}

%TCIMACRO{\TeXButton{Proof}{\proof}}%
%BeginExpansion
\proof%
%EndExpansion
This follows from Theorems \ref{compression theorem} and \ref{extend Jones
theorem}.%
%TCIMACRO{\TeXButton{End Proof}{\endproof}}%
%BeginExpansion
\endproof%
%EndExpansion

\begin{figure}[th]
\begin{center}
\includegraphics[width=5.1in,height=0.53in]{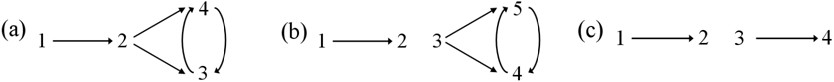}
\end{center}
\caption{We use (b) and (c) to find a $G$-sourcing set for the relation of (a).}
\label{Two-Step Figure}
\end{figure}

The directed graphs in Figure \ref{Two-Step Figure}\ represent reflexive
relations. The vertices have been replaced with elements of the sets they
represent.

\begin{example}
\label{stable to partial order example}Suppose $\rho _{1}$ is the stable
relation on $X_{1}=\left\{ 1,2,3,4\right\} $ determined by (a) in Figure \ref%
{Two-Step Figure}. The only clasp in $X_{1}$ is 2, which is an unlocked
clasp, so $X_{1}$ is the compression of a preordered set by part 1 of
Theorem \ref{Stable to Transitive Thm}. If $\rho _{2}$ is the reflexive
relation on $X_{2}=\left\{ 1,2,3,4,5\right\} $ determined by (b) in Figure %
\ref{Two-Step Figure} then $\rho _{2}$ is a preorder on $X_{2}$ and there is
a compression map $\theta :X_{2}\rightarrow X_{1}$ given by $\theta \left(
1\right) =1$, $\theta \left( 2\right) =2$, $\theta \left( 3\right) =2$, $%
\theta \left( 4\right) =3$, and $\theta \left( 5\right) =4$. By equating
paired elements we may relate $\left( X_{2},\rho _{2}\right) $ to a partial
order. This gives $\rho _{3}$, the minimally-connected partial order on $%
X_{3}=\left\{ 1,2,3,4\right\} $ determined by (c) in Figure \ref{Two-Step
Figure}. By Theorem \ref{Hasse Diagram Lemma} a $G$-grading set for $\rho
_{3}$ is $\sigma _{3}=\left\{ \left( 1,2\right) ,\left( 3,4\right) \right\} $%
. Following the proof of Theorem \ref{extend Jones theorem} a $G$-grading
set for $\rho _{2}$ is $\sigma _{2}=\beta \cup \gamma $ with $\beta =\left\{
\left( 1,2\right) ,\left( 3,4\right) \right\} $ and $\gamma =\left\{ \left(
4,5\right) \right\} $. The proof of Theorem \ref{compression theorem} shows
a $G$-grading set of $\rho _{1}$ is $\sigma _{1}=\left\{ \left( 1,2\right)
,\left( 2,3\right) ,\left( 3,4\right) \right\} $. Thus for any associative
ring with unity $R$ and any group $G$ the $G$-gradings of $S=I\left( X,\rho
,R\right) $ are uniquely determined by functions from $\sigma _{1}$ to $G$.
\end{example}

Suppose $\sigma $ is the Hasse diagram for a partial order $\left( X,\leq
\right) $ and $F$ is a field. The incidence algebra $A=I\left( X,\leq
,F\right) $ is said to have the free-extension property (see \cite[%
Definition 2]{MJ}) if there exists $\gamma \subseteq \sigma $ such that for
any group $G$ and any function $\phi :\gamma \rightarrow G$, there is a
unique grading of $I\left( X,\leq ,F\right) $ such that $e_{xy}\in A_{\phi
\left( x,y\right) }$ for all $x,y\in X$ such that $\left( x,y\right) \in
\gamma $.

Theorem \ref{All Good Theorem} shows a good group grading of a generalized
incidence ring induced by a homomorphism must be a finite grading. This kind
of result was already known for incidence algebras over partial orders (see 
\cite[Theorem 3.3]{MJD}). The free-extension property may fail for partial
orders on an infinite set even if it contains a cross-cut of length one. We
describe a partial order with a minimal element which does not have the
free-extension property.

\begin{example}
\label{infinite support example} Suppose $F$ is a field, $G$ is a group, and 
$X$ is the set of natural numbers with $\rho $ the usual ordering. The set $%
\sigma =\left\{ \left( m,m+1\right) :m\in X\right\} $ is the arrow set of
the Hasse diagram of $X$ hence $\sigma $ is a $G$-grading set for $X$ by
Theorem \ref{Hasse Diagram Lemma}. Suppose there is an element $g\in G$ of
infinite order. We let $\phi :\sigma \rightarrow G$ be the function given by 
$\phi \left( m,m+1\right) =g^{m}$ for all $m\in X$. Then there is a unique
homomorphism $\Phi :\rho \rightarrow G$ such that $\Phi |_{\sigma }=\phi $,
but this does not determine a grading for $I\left( X,\rho ,F\right) $ since $%
\func{Im}\Phi $ is infinite. Therefore $I\left( X,\rho ,F\right) $ does not
have the free-extension property even though it satisfies the requirements
of \cite[Theorem 4]{MJ}.
\end{example}

\paragraph{Acknowledgements}

A sabbatical appointment supported a visit to the Atlantic Algebra Center of
Memorial University in Newfoundland Canada. It is a pleasure to thank the
faculty for their hospitality. This paper answers only some of the many fine
questions they posed to the author. The author also thanks Martin Erickson
of Truman State University and Gene Abrams of the University of Colorado, at
Colorado Springs for suggesting improvements to the writing in earlier
drafts.

\end{document}